\documentclass[12pt,leqno]{article}
\usepackage{amsmath,amsthm}
\usepackage{amsfonts}
\newtheorem{theorem}{Theorem}
\newtheorem{lemma}[theorem]{Lemma}

\newcommand{\half}{{\mbox{\small $\frac{1}{2}$}}}
\begin{document}
\bibliographystyle{plain}
\begin{center}
{ \LARGE  Splitting: Tanaka's SDE revisited}
\\ ~ \\
 { \large \sl J. WARREN\footnote{University of Warwick, United Kingdom}}
\end{center}

What follows is my attempt to understand a set of ideas being
developed by Boris Tsirelson. I do this by studying a
specific, and I hope interesting, example.

 Tanaka's SDE
is one of the easiest examples of a stochastic differential equation with
no strong solution. Suppose
$\bigl( X_t; t
\geq 0 \bigr)$ is a real-valued Brownian motion starting from zero and 
we put
$B_t=\int_0^t  sgn(X_s) dX_s$ then $B$ is also a Brownian motion  and
Tanaka's SDE
\begin{equation}
\label{tan}
X_t=  \int_0^t  sgn(X_s) dB_s,
\end{equation}
is satisfied. But the trajectory of $B$ does not determine that of $X$.
Recall Tanaka's formula 
\begin{equation}
\label{tan2}
\vert X_t \vert= \int_0^t sgn(X_s) dX_s + L_t,
\end{equation}
where $\bigl(L_t;t\geq 0\bigr)$ is the local time process of $X$ at zero.
 We find 
\begin{equation}
\label{tan3}
  \vert X_t \vert=B_t + \sup_{s\leq t} (-B_s) 
\end{equation}
but $B$ does not tell us the signs of the excursions from zero made by
$X$.
 
In a discrete time framework things, work out differently. Let $\bigl( X_n;
n \geq 0 \bigr)$ be the symmetric nearest neighbour random walk on the
integers. Define $sgn(a)$ to be $+1$ if $a \geq 0$ and $-1$ if $a<0$,
and then let $Z_0=0$ and $ Z_{n+1}-Z_n= sgn(X_n) \bigl(
X_{n+1}-X_n
\bigr)$ then $\bigl(Z_n; n\geq 0)$ is again a symmetric random walk and 
we may write a discrete version of equation \eqref{tan}:
\begin{equation}
\label{tand}
X_n = \sum_{k=0}^{n-1} sgn(X_k) \bigl(
Z_{k+1}-Z_k
\bigr) .
\end{equation} 
The equations \eqref{tan2} and \eqref{tan3} have discrete time versions:
\begin{equation}
\vert X_n+\half \vert -\half= \sum_{k=0}^{n-1} sgn(X_k) \bigl(
X_{k+1}-X_k\bigr) +  L_n
\end{equation} 
where $L_0=0$ and for $n \geq 1$ we define $L_n =\sum_{k=0}^{n-1} 1_{\bigl(
X_k, X_{k+1}\in\{0,-1\}\bigr)}$,
and,
\begin{equation}
\label{tand3}
\vert{X_n}+ \half \vert - \half=  Z_n + \sup_{k\leq n} \bigl( -Z_k \bigr).
\end{equation}
The halves appear because of the lack of symmetry in our definition of
$sgn$- it is not something to worry about.  This time $Z$ does determine
$X$: the information about whether
$X_n$ is below or above $-\half$- which at first sight appears to be missing 
in \eqref{tand3}- is coded in the following way. Find the last time
$r \in \{0, 1, \ldots ,n\} $ such that $Z_r= -\sup_{k\leq r} \bigl( -Z_k
\bigr)$. If this $r$ is even then $X_n+\half$ is positive while if it is
odd then
$X_n+\half$ is negative.

This note is concerned with understanding what happens to this precious
information about the sign of $X$ when we try to obtain the
continuous-time model by taking scaling limits of the discrete model. 
This is  inspired by  work of Boris Tsirelson on the spectra of
noises and  stability- see \cite{tsir3}, \cite{tsir2}, and \cite{sch}. 

One way to understand why the information about the signs does not
survive the passage to the limit is to observe that it is noise
sensitive. Instead of one copy of the random walk take a pair
$(Z, Z^\prime)$ that are $\rho$-correlated where $\rho \in (0,1)$. This
means that together they form
a nearest-neigbour random walk on the lattice ${\mathbb Z}^2$ with
${\mathbb E}\bigl[ (Z_{k+1}-Z_{k})(Z^\prime_{k+1}-Z^\prime_{k})
\bigr]=\rho$.  We think of $Z^\prime$ as being a perturbation of $Z$.
Now construct $X$ and $X^\prime$ from $Z$ and $Z^\prime$ so that
equation \eqref{tand} and its prime version hold. 
As $n$ becomes large  ( as it will when we try to take scaling limits)
we find that $sgn(X_n)$ and $sgn(X^\prime_n)$ become uncorrelated no
matter how strong the (fixed) correlation $\rho$ between $Z$ and
$Z^\prime$ is.  Thus, in a certain sense, $sgn(X_n)$ is
asymptotically sensitive to  small perturbations of $Z$.

The discussion of the previous paragraph, although very elementary, is (a
variant of) the observation that eventually led Tsirelson to profound
results in the theory of filtrations \cite{tsir4}. There good account
of the story in \cite{sch}.

Next instead of constructing $Z^\prime$ by perturbing the
whole path of $Z$ we may only perturb some sections. More precisely
let $A \subset [0,1]$ be a finite union of closed intervals with (to be
safe) dyadic rational end-points. Fix $\rho \in (0,1)$. For each $n$
construct a random walk $\bigl((Z_k, Z^\prime_k); 0 \leq k \leq n \bigr)$
on 
${\mathbb Z}^2$ with
\begin{equation*}
{\mathbb E}\bigl[ (Z_{k+1}-Z_{k})(Z^\prime_{k+1}-Z^\prime_{k}) \bigr ]=
\begin{cases}
\rho & k2^{-n} \in A \\
1 &  k2^{-n} \in A^c.
\end{cases}
\end{equation*}
Now as before construct $X$ and $X^\prime$ and consider the correlation
of $sgn(X_n)$ and $sgn(X^\prime_n)$. This time as $n$ tends towards
infinity we obtain a nontrivial limit which we denote by $\phi(\rho,A)$. 
  
We write the Wiener chaos expansion of any
random variable belonging to ${\cal L}^2(B)$ in the form
\[
\hat{f}_0 + \int_0^1  \hat{f}_1(s) dB_s +  \int_0^1 \int_0^{s_2}
\hat{f}_2(s_1,s_2) dB_{s_1}dB_{s_2} + \ldots
\]
Then we construct a finite measure on $\cup_{n\geq 0} \bigl \{ (s_1, s_2 , 
\ldots ,s_n) \in [0,1]^n \vert s_1 <s_2 \ldots <s_n \bigr\}$
having density  $\vert \hat{f}_n(s_1,\ldots,s_n) \vert ^2$ with respect to 
Lebesgue
measure. We call this the spectral measure of the random
variable whose chaos expansion we used. If we start with variable having
 ${\cal L}^2$-norm equal to one this measure is a probability measure-
and we can think of it as determining the law of a finite random subset
$S$ of $[0,1]$. Thus there is a probability $\vert\hat{f}_0\vert^2$ that $S$ is
empty, a probability $\vert\hat{f}_1(s)\vert^2ds$ that it contains a single 
point
lying in $(s,s+ds)$ and so on. Suppose for a moment that the information on 
signs did survive in the
limit, and that the $X$ satisfying Tanaka's SDE was some functional of
the Brownian motion $B$. Then we could apply this construction to
$sgn(X_1)$- and obtain a random subset $S$.
If $A$ is a fixed subset of $[0,1]$ once more then let $\vert S \cap A
\vert$ denote the number of  points of $S$ lying in $
A$. Then it is reasonable to expect that  $\phi(\rho,A)$ would  be given by
\begin{equation}
\label{rep}
\phi(\rho,A)= {\mathbb E} \left[ \rho^{\vert S \cap A \vert} \right ].
\end{equation}
The ${\mathbb E}$ appearing here is with respect to the law of $S$ which
does not live on the
same probability space as $X$ and $B$.

The discussion of the preceding paragraph is based on a false
premise, but nevertheless there is a random subset-
still denoted by
$S$- such that equation
\eqref{rep} holds. This subset possess, with probability one, 
an infinite number of elements. If we take $A= [0,1]$ then $\vert S
\cap A \vert$ is infinity and $\rho^\infty=0$ (by definition if you
like!). Thus $\phi(\rho,[0,1])$ is $0$ for any $\rho$- this is just the 
sensitivity to noise property with which we began. In what follows we
examine the law of this $S$ more closely. Not surprisingly the Wiener chaos
expansion is our principle tool.

For any $x >0$ and $t \in (0,1)$ let $m_{(t,x)}$ denote the spectral
measure  of
\[
1\left( \sup_{h\in [t,1]}(B_t-B_h) <x\right).
\]
Note that  the total mass of this measure is just ${\mathbb P}\left(
\sup_{u \in [t,1]}(B_t-B_u) <x\right)<1$, but we will nevertheless speak
of a random set $S$ having distribution  $m_{(t,x)}$. This subset is
supported on $[t,1]$.  Let $q_h(x,dy)$ denote the (defective) transition
probability distributions of Brownian motion killed on hitting $0$. 

\begin{lemma} Suppose that $0<s<t<1$ and that $S$ is distributed
according to $m_{(s,x)}$. Then the subset $S\cap[t,1]$ is distributed
according to 
\[
\int q_{t-s}(x,dy) m_{(t,y)}.
\]
\end{lemma}
Now suppose that $\bigl(\lambda_t(dx); 0<t\leq 1\bigr)$ is an entrance
law for killed Brownian motion; thus, for any $0<s<t\leq 1$,
\[
\lambda_t(dy)=\int q_{t-s}(x,dy) \lambda_s(dx).
\]
We may define a family of measures $m_{(t,\lambda)}$ for $t\in (0,1)$ via
\[
m_{(t,\lambda)}= \int \lambda_t(dy) m_{(t,y)},
\]
and by virtue of the lemma they have a certain consistency property-
that is - if $S$ is distributed according to $m_{(s,\lambda)}$ then for
any $t>s$ the intersection $S
\cap [t,1]$ is distributed according to $m_{(t,\lambda)}$. Note that the
total mass of $m_{(t,\lambda)}$ is the `probability' that the killed BM
survives to time $1$ when it is started according to $\lambda$- this does
not depend on $t$.  Because of this consistency property there is a
random set $S\subset (0,1]$ whose distribution we denote by $m_\lambda$ 
whose intersection with $[t,1]$ has distribution  $m_{(t,\lambda)}$ for
any $t$. This $S$ may  be infinite- there is the possibility of $0$
being an accumulation point.

 From this point on we will take $\bigl(\lambda_t; 0<t\leq 1\bigr)$ to be
a multiple of the entrance law for the It\^{o} excursion measure for the
positive excursions of Brownian motion. Choose this multiple so that
$m_\lambda$ becomes a probability measure. More explicitly we have
\begin{equation}
\lambda_t(dy)= y\;t^{-3/2}\exp\left\{- \frac{y^2}{2t}\right\} \;
dy \qquad \qquad y>0.
\end{equation}
Recall that a random variable is said to be arc-sine distributed if it
has distribution
\[
s(dt)=\frac{dt}{\pi \sqrt{t(1-t)}}1_{[0,1]}(t)\, dt.
\] 
The time at which a  BM  attains its
minimum between times $0$ an $1$ is so distributed. 
\begin{theorem}
The limits $\phi(\rho,A)$ exist and admit the following description.
Take two random subsets $S_1$ and $S_2$ distributed according to 
$m_\lambda$ and a $[0,1]$-valued random variable $g$ with the arc-sine
distribution . Suppose that $S_1$, $S_2$ and $g$ are independent.  Take
\[
S= g\bigl(1-S_1\bigr) \cup \bigl( (1-g)S_2 +g\bigr),
\]
then for all $\rho$ and $A$
\[
\phi(\rho,A)= {\mathbb E} \left[ \rho^{\vert S \cap A \vert} \right ].
\]
\end{theorem}
\begin{proof}[Proof of Lemma]
Begin by writing
\begin{multline*}
1\left( \sup_{h \in [s,1]}(B_s-B_h) <x\right) =
1\left( \sup_{h \in [s,t]}(B_s-B_h) <x\right) \\
\times 1\left( \sup_{h \in
[t,1]}(B_t-B_h) <x+B_t-B_s\right).
\end{multline*}
Condition on $\bigl( B_r; r \leq t \bigr)$ and then replace the second factor  
with its Wiener chaos expansion and so
obtain an expansion of which the typical term is
\begin{multline*}
\int_t^1dB_{h_1}\int_t^{h_1} dB_{h_2} \ldots \int_t^{h_{k-1}}dB_{h_k}
\; \\
 1\left( \sup_{h \in [s,t]}(B_s-B_h) <x\right)
\hat{f}_k(t,x+B_t-B_s\vert h_1,\ldots ,h_k).
\end{multline*}
We now replace each integrand  by its chaos expansion- this must simply 
result in the chaos expansion of
\[
1\left( \sup_{h \in [s,1]}(B_s-B_h) <x\right). 
\]
On comparing the two expansions it may be seen  that if $S$ is distributed
according to $m_{(s,x)}$ then $S \cap [t,1]$ contains exactly k points at 
positions
$(h_1,h_1+dh_1)$ through $(h_k,h_k+dh_k)$ with probability
\[
{\mathbb E} \left[ 1\left( \sup_{h \in [s,t]}(B_s-B_h) <x\right)
\vert\hat{f}_k(t,x+B_t-B_s\vert h_1,\ldots ,h_k)\vert^2 \right]dh_1
\ldots dh_k,
\]
but since $ \vert\hat{f}_k(t,y\vert h_1,\ldots ,h_k)\vert^2dh_1
\ldots dh_k$ is just the probability distribution of $S$ under 
$m_{(t,y)}$ we are done.
\end{proof}
\begin{proof}[Proof of Theorem]
{\em  Stage 1.} Fix an admissible subset $A$.   For each $n$ consider the
correlated random walk $\bigl((Z_k, Z^\prime_k); 0 \leq k \leq n \bigr)$.
There is the usual weak convergence in the space of continuous ${\mathbb
R}^2$-valued paths to a process $\bigl( (B_t,B^\prime_t); 0\leq t \leq
1\bigr)$, each component of which forms a one-dimensional Brownian motion
and their co-variation is simply:
\begin{equation*}
dB_tdB^\prime_t=
\begin{cases}
\rho dt & t \in A \\
dt &  t \in A^c.
\end{cases}
\end{equation*}
 Let $g$ be the time at which $B$ attains its minimum between times $0$
and $1$, and similarly define $g^\prime$. Now the correlation of
$sgn(X_n)$ and
$sgn(X^\prime_n)$ can be split into the sum of two contributions. One
arises when  the random walks $Z$ and $Z^\prime$ attain their minimum
(between times $0$ and $n$ ) values simultaneously - in this case
$sgn(X_n)$ and $sgn(X^\prime_n)$ are equal.
The remaining contribution  tends to zero  for large $n$   and so 
 the limits $\phi(\rho,A)$ exist and are given by
\[
\phi(\rho,A)= {\mathbb E} \left[ 1_{ (g=g^\prime)} \right].
\]
{\em Stage 2.} The two random times $g$ and $g^\prime$ can only be equal
if their common value lies in one of the components of $A^c$. For each
such component we consider the common time  at which $B$ and
$B^\prime$ attain their respective minimum (over that component) and
compute the probability that this is actually the global minimum of both
Brownian motions. We obtain
\begin{equation*}
{\mathbb E} \left[ 1_{ (g=g^\prime)} \right]=\int_{A^c}  \int_0^\infty
 \int_0^\infty  p(u_t,v_t;dt,dy_1,dy_2)  m_{(v_t,y_2)}
\bigl[ \rho^{\vert S \cap A\vert } \bigr] m_{(1-u_t,y_1)}
\bigl[ \rho^{\vert (1-S) \cap A\vert } \bigr],
\end{equation*}
where
\begin{align*}
u_t &= \sup \{ h<t: h \in A\},
\\
v_t &= \inf \{ h>t : h \in A\},
\end{align*}
and 
$ p(u,v;dt,dy_1,dy_2)$ is the law of the triple

\[
 \bigl( g(u,v), B_u-I(u,v), B_v-I(u,v) \bigr),
\]
 $g(u,v)$ denoting
the time at which $B$ attains its minimum $I(u,v)= \inf \{ h \in[u,v]:
B_h\}.$

{\em Stage 3.} 
By virtue of the scaling properties of BM  we have
\begin{align*}
  m_{(v,y)}
\bigl[ \rho^{\vert S \cap A\vert } \bigr]  
&= m_{\left((v-t)/(1-t),y/\sqrt{1-t}\right)}
\bigl[ \rho^{\vert ( t+(1-t)S )\cap A\vert } \bigr]
\\
m_{(1-u,y)}
\bigl[ \rho^{\vert (1-S) \cap A\vert } \bigr] &=
m_{\left((t-u)/t,y/\sqrt{t}\right)}
\bigl[ \rho^{\vert t(1-S) \cap A\vert } \bigr].
\end{align*}
A well-known exercise (Revuz and Yor \cite{ry}, chapter XII) 
confirms that
\begin{align*}
p(u,v;dt,dy_1,dy_2) &= \frac {dt}{\pi}
\lambda_{t-u}(dy_1)\lambda_{v-t}(dy_2) 
\\
&= s(dt)
\lambda_{(t-u)/t}(dy_1/\sqrt{t})\lambda_{(v-t)/(1-t)}(dy_2/\sqrt{1-t}).
\end{align*}
Putting these into the formula obtained in the previous section and
recalling the definition of $m_\lambda$  we
obtain the desired result:
\[
{\mathbb E} \left[ 1_{ (g=g^\prime)} \right]=\int_{A^c}  
s(dt) m_\lambda \bigl[ \rho^{\vert ( t+(1-t)S )\cap A\vert } \bigr] 
m_\lambda \bigl[ \rho^{\vert t(1-S) \cap A\vert } \bigr].
\]
\end{proof}

It is possible to generalise the family of measures $m_{(t,x)}$ from which
we obtained $m_\lambda$. Starting from a bounded function $f$ defined on
${\mathbb R}_+$ we may  expand 
\[
f\bigl(B_1-B_t+x\bigr) 1\left( \sup_{h\in [t,1]}(B_t-B_h) <x\right),
\]
and whence construct a measure $m^f_\lambda$. I would like to know 
when such measures corresponding to different $f$ are equivalent and in
this case how to compute the Radon-Nikod\'{y}m density. 
This is part of the problem of obtaining the spectral resolution (see
\cite{tsir1}) of the  noise of splitting. This is a noise richer than white 
noise: in addition to the increments of a Brownian motion $B$ it carries a 
countable collection of independent Bernoulli random variables which are 
attached to the local minima of $B$.

\end{document}